\documentstyle[12pt]{article}
\newcommand\be{\begin{equation}}
\newcommand\ee{\end{equation}}
\newcommand\bea{\begin{eqnarray}}
\newcommand\eea{\end{eqnarray}}
\newcommand\ket[1]{|#1\rangle}

\newcommand\braket[2]{\langle #1|#2\rangle}
\newcommand{\fatalpha}{{\bf \alpha \kern -0.44em \alpha}}
\newcommand{\fatsigma}{{\bf \sigma \kern -0.54em \sigma}}
\newcommand{\tpchi}{{\bf \chi \kern -0.35em \chi}}
\newcommand{\llambda}{{\bf \lambda \kern -0.45em \lambda}}

\renewcommand{\theequation}{\arabic{equation}}
\renewcommand{\theequation}{\thesection-\arabic{equation}}
\bibliography{plain}
\title{\bf Recursive calculation of effective resistances in distance-regular networks based on Bose-Mesner algebra and Christoffel-Darboux identity}\vspace{20mm}
\author{ M. A. Jafarizadeh$^{a,b,c}$
 \thanks{E-mail:jafarizadeh@tabrizu.ac.ir},
 R. Sufiani$^{a,c}$
 \thanks{E-mail:sofiani@tabrizu.ac.ir},
 S. Jafarizadeh$^{d}$
\\ $^a${\small Department of Theoretical Physics and Astrophysics,
The University of Tabriz, Tabriz 51664, Iran.} \\ $^b${\small
Institute for Studies in Theoretical Physics and Mathematics, Tehran
19395-1795, Iran.} \\ $^c${\small Research Institute for Fundamental
Sciences, Tabriz 51664, Iran.}\\ $^d${\small Department of
Electrical and computer engineering, The University of Tabriz,
Tabriz 51664, Iran.}} \pagebreak

\vspace{20mm}
\begin{document}
\maketitle \vspace{15mm}
\newpage
\begin{abstract}
Recently in \cite{jss1}, the authors have given a method for
calculation of the effective resistance (resistance distance) on
distance-regular networks, where the calculation was based on
stratification introduced in \cite{js} and Stieltjes transform of
the spectral distribution (Stieltjes function) associated with the
network. Also, in Ref. \cite{jss1} it has been shown that the
resistance distances between a node $\alpha$ and all nodes $\beta$
belonging to the same stratum with respect to the $\alpha$
($R_{\alpha\beta^{(i)}}$, $\beta$ belonging to the $i$-th stratum
with respect to the $\alpha$) are the same. In this work, an
algorithm for recursive calculation of the resistance distances in
an arbitrary distance-regular resistor network is provided, where
the derivation of the algorithm is based on the Bose-Mesner
algebra, stratification of the network, spectral techniques and
Christoffel-Darboux identity. It is shown that the effective
resistance on a distance-regular network is an strictly increasing
function of the shortest path distance defined on the network. In
the other words, the two-point resistance
$R_{\alpha\beta^{(m+1)}}$ is strictly larger than
$R_{\alpha\beta^{(m)}}$. The link between the resistance distance
and random walks on distance-regular networks is discussed, where
the average commute time (CT) and its square root (called
Euclidean commute time (ECT)) as a distance are related to the
effective resistance. Finally, for some important examples of
finite distance- regular networks, the resistance distances are
calculated.

{\bf Keywords: resistance distance, association scheme,
stratification, distance-regular networks, Christoffel-Darboux
identity}

{\bf PACs Index: 01.55.+b, 02.10.Yn }
\end{abstract}

\vspace{70mm}
\newpage
\section{Introduction}
A classic problem in electric circuit theory studied by numerous
authors over many years is the computation of the resistance
between two nodes in a resistor network (see, e.g.,
\cite{Cserti}). The study of electric networks was formulated by
Kirchhoff \cite{8} more than 150 years ago as an instance of a
linear analysis. Besides being a central problem in electric
circuit theory, the computation of resistances is also relevant to
a wide range of problems ranging from random walks (see
\cite{Doyle}), the theory of harmonic functions \cite{4}, to
lattice Green's functions \cite{6,6a,6b,6c,6d}. The connection
with these problems originates from the fact that electrical
potentials on a grid are governed by the same difference equations
as those occurring in the other problems. For this reason, the
resistance problem is often studied from the point of view of
solving the difference equations, which is most conveniently
carried out for infinite networks. In the case of Green's function
approach, for example, past efforts \cite{Cserti}, \cite{77} have
been focused mainly on infinite lattices. Little attention has
been paid to finite networks, even though the latter are those
occurring in real life.

Within the theory of electrical networks, the standard method to
compute the two-point resistances on a network is via the
Moore-Penrose generalized inverse \cite{moore} or pseudo inverse
of the Laplacian $L$ of the underlying network, where the
Laplacian is a matrix whose off-diagonal entries are the
conductances connecting pairs of nodes. Just as in graph theory
where everything about a graph is described by its adjacency
matrix (whose elements is 1 if two vertices are connected and 0
otherwise), everything about an electric network is described by
its Laplacian. It should be noticed that, the concept of two-point
resistance called also effective resistance and resistance
distance and molecular structure-descriptors based on them, were
much studied in the chemical literature \cite{12'}-\cite{23'}.
Also, it is closely related to the average first-passage time and
the average commute time which are two important quantities in
random walk models defined based on the Markov chain. It is shown
in \cite{fouss} that the computation of the average commute time
can be obtained from the pseudo inverse of the Laplacian $L$
($L^{-1}$) of the underlying network. Also, it has been shown that
this quantity and its square root is a distance, since it can be
shown \cite{fouss} that $L^{-1}$ is symmetric and positive
semidefinite. It is therefore called the Euclidean Commute Time
(ECT) distance. In fact the ECT distance is the same as resistance
distance or effective resistance (the effective resistance is
symmetric and satisfies the triangle inequality and so is a
distance metric). Therefore, any clustering algorithm
(hierarchical clustering, $k$-means, etc) which can be used in
conjunction with the ECT distance, deals with the effective
resistance. Also, Laplacian eigenmaps which is one of the
graph-based approaches on dimensionality reduction and manifold
learning and recently proposed by Belkin and Niyogi in
\cite{belkin} and share many similarities with other recent
spectral algorithms for clustering and embedding of data, for
instance Kernel PCA (Principle Component Analysis) algorithm
\cite{smola} and spectral methods for image segmentation
\cite{meila} (for a unifying view of the behavior of spectral
embeddings and clustering algorithms, see \cite{brand}) deals with
the Laplacian of the graph assigned to the data and uses the
resistance distance. One of the most important aspects of spectral
methods for clustering and embedding, including Laplacian
eigenmaps, is the fact that they are all posed as eigenvalue
problems. But clearly, for too large matrices, the computation by
pseudo inverse becomes cumbersome.

Recently, the authors have given a method for calculation of the
resistance distance on distance-regular networks \cite{jss1},
where the calculation is based on stratification introduced in
\cite{js} and Stieltjes transform of the spectral distribution
(Stieltjes function) associated with the network. Also, in
Ref.\cite{jss1} it has been shown that the resistances between a
node $\alpha$ and all nodes $\beta$ belonging to the same stratum
with respect to the $\alpha$ ($R_{\alpha\beta^{(i)}}$, $\beta$
belonging to the $i$-th stratum with respect to the $\alpha$) are
the same and the analytical formulas have been given for two-point
resistances $R_{{\alpha\beta^{(i)}}},i=1,2,3$ in terms of the size
of the network and corresponding intersection array without any
need to know the spectrum of the pseudo inverse $L^{-1}$. In this
work, our starting point is along the same line by considering the
Laplacian matrix associated with a network, but we use the
algebraic structure of distance-regular networks (Bose-Mesner
algebra) such as stratification and spectral techniques specially
the well known Christoffel-Darboux identity \cite{tsc} from the
theory of orthogonal polynomials to give a recursive formula for
calculation of all of the resistance distances
$R_{{\alpha\beta^{(i)}}}$, $i=1,2,...,d$ ($d$ is the diameter of
the graph) on the network without any need to the spectrum of the
pseudo inverse $L^{-1}$. In fact, in order to evaluate the
resistance distance, one needs only to know the intersection array
of the network. The other main result of the derived recursive
formula is that, the resistance distance strictly increases by
increasing the shortest path distance defined on the network,
i.e., $R_{{\alpha\beta^{(m+1)}}}-R_{{\alpha\beta^{(m)}}}>0$ for
all $m=1,2,...,d-1$.

The organization of the paper is as follows. In section $2$, we give
some preliminaries such as association schemes, distance-regular
networks, stratification of these networks and Stieltjes function
associated with the network. Section $3$ is devoted to the concept
of two-point resistance on regular networks and its relation with
random walks. In Section $4$ (main section of the paper) we show
that the two-point resistance increases with increasing of the
number of stratum and give a recursive formula for calculation of
two-point resistances $R_{\alpha\beta^{(i)}}$ for $i=1,2,...d$ in
distance-regular networks, where the algorithm is based on the
Bose-Mesner algebra and spectral techniques specially the
Christoffel-Darboux identity. Section $5$ is devoted to calculation
of effective resistances on some examples of distance-regular
networks. The paper is ended with a brief conclusion and an appendix
containing calculation of two-point resistances
$R_{\alpha\beta^{(i)}}$, $i=1,2,...,d$ for some important finite
distance-regular networks.
\section{Preliminaries}
In this section we give some preliminaries such as definitions
related to association schemes, corresponding stratification,
distance-regular networks and Stieltjes function associated with a
distance-regular network.
\subsection{Association schemes}
First we recall the definition of association schemes. The reader
is referred to Ref.\cite{Ass.sch.}, for further information on
association
schemes.\\
\textbf{Definition 2.1} (Symmetric association schemes). Let $V$
be a set of vertices, and let $R_i(i = 0, 1,..., d)$ be nonempty
relations on $V$ (i.e., subset of $V\times V$). Let the following
conditions (1), (2), (3) and (4) be satisfied. Then, the relations
$\{R_i\}_{0\leq i\leq d}$
on $V\times V$ satisfying the following conditions\\
$(1)\;\ \{R_i\}_{0\leq i\leq d}$ is a partition of $V\times V$\\
$(2)\;\ R_0=\{(\alpha, \alpha) : \alpha\in V \}$\\
$(3)\;\ R_i=R_i^t$ for $0\leq i\leq d$, where
$R_i^t=\{(\beta,\alpha) :(\alpha, \beta)\in R_i\} $\\
$(4)$ For $(\alpha, \beta)\in R_k$, the number  $p^k_{ij}=\mid
\{\gamma\in V : (\alpha, \gamma)\in R_i \;\ and \;\
(\gamma,\beta)\in R_j\}\mid$ does not depend on $(\alpha, \beta)$
but only on $i,j$ and $k$,\\ define a symmetric association scheme
of class $d$ on $V$ which is denoted by $Y=(V,\{R_i\}_{0\leq i\leq
d})$. Furthermore, if we have $p^k_{ij}=p^k_{ji}$ for all
$i,j,k=0,2,...,d$, then $Y$ is called commutative.

The number $v$ of the vertices, $|V|$, is called the order of the
association scheme and $R_i$ is called $i$-th relation (colour).
For example, in the resistor networks the relations $R_i$,
$i=0,1,...,d$ can be interpreted as $d+1$ different kinds of
conductances, i.e., two nodes $\alpha,\beta$ have the $i$-th
relation with each other if and only if the conductance between
them be $c_i$ (see Figure 1). In this paper we will deal with the
special case where the conductance between two nodes
$\alpha,\beta$ is $c\equiv c_1$ if the nodes be adjacent, i.e.,
$(\alpha,\beta)\in R_1$ and the other conductances corresponding
to the other relations will be taken to zero.

Also note that, the intersection number $p_{ij}^k$ can be
interpreted as the number of vertices which have relation $i$ and
$j$ with vertices $\alpha$ and $\beta$, respectively provided that
$(\alpha,\beta)\in R_k$, and it is the same for all elements of
relation $R_k$. For all integers $i$ ($0 \leq i \leq d$), set
$\kappa_i=p_{ii}^{0}$ and note that $\kappa_i\neq 0$, since $R_i$
is non-empty. We refer to $\kappa_i$ as the $i$-th valency of $Y$.

 For examples of association schemes, consider a cube known
as Hamming scheme $H(3,2)$, in which $V$ (the vertex set) is the
set of $3$-tuples with entries in $F_2=\{0,1\}$. Two vertices are
connected if and only if they differ by exactly one entry (see
Figure 1). The distance between vertices, i.e. the length of the
shortest edge path connecting them, will then indicate which
relation they are contained in. E.g., if $x=(0,0,1)$, $y=(0,1,1)$
and $z=(1,0,1)$, then $(x,y)\in R_1$, $(x,z)\in R_1$ and $(y,z)\in
R_2$. As an another example, consider the octahedron (a special
case of a square dipyramid with equal edge lengths) which is the
same as Johnson scheme $J(4,2)$, in which the vertex set $V$
contains all $2$-element subsets of the set $\{1,2,3,4\}$ and two
vertices are adjacent if and only if they intersect in exactly one
element. Two vertices are then at distance $i$, $i=0,1,2$ if and
only if they have exactly $2-i$ elements in common (see Figure 2).

Let $Y=(X,\{R_i\}_{0\leq i\leq d})$ be a commutative symmetric
association scheme of class $d$, then the matrices
$A_0,A_1,...,A_d$ defined by
\begin{equation}\label{adj.}
\bigl(A_{i})_{\alpha, \beta}\;=\;\cases{1 & if $\;(\alpha,
\beta)\in R_i$ \cr 0 & otherwise\cr},
\end{equation}
are adjacency matrices of $Y$ such that
\begin{equation}\label{ss}
A_iA_j=\sum_{k=0}^{d}p_{ij}^kA_{k}.
\end{equation}
From (\ref{ss}), it is seen that the adjacency matrices $A_0, A_1,
..., A_d$ form a basis for a commutative algebra \textsf{A} known
as the Bose-Mesner algebra of $Y$. Since the matrices $A_i$
commute, they can be diagonalized simultaneously. The Bose-Mesner
algebra has a second basis $E_0, ...,E_d$, such that, $E_iE_j=
\delta_{ij}Ei$ and $\sum_{i=0}^d E_i=I$ with $E_0=1/n J$ ($J$ is
the all-one matrix) \cite{Ass.sch.}. The matrices $E_i$, for
$(0\leq i, j\leq d)$ are known as the primitive idempotents of the
$Y$. Then, there are matrices $P$ and $Q$ such that the two bases
of the Bose-Mesner algebras can be related to each other as
follows
$$A_i=\sum_{j=0}^d{\mathrm{P}}_{ij}E_j,\;\;\ 0\leq j\leq d,$$
\begin{equation}\label{ap}
E_i=\frac{1}{N}\sum_{j=0}^d{\mathrm{Q}}_{ij}A_j,\;\;\ 0\leq j\leq
d.
\end{equation}
where, $N$ denotes the cardinality of $X$. Then clearly
\begin{equation}\label{pq}
\mathrm{P}\mathrm{Q}=\mathrm{Q}\mathrm{P}=NI.
\end{equation}
It also follows that
\begin{equation}\label{eign}
A_jE_i={\mathrm{P}}_{ij}E_i,
\end{equation}
which indicates that the ${\mathrm{P}}_{ij}$ is the $i$-th
eigenvalue of $A_j$ and that the columns of $E_i$ are the
corresponding eigenvectors. Also, $m_i:=tr E_i=N\langle
\alpha|E_i|\alpha\rangle$ is the rank of the idempotent $E_i$
which gives the multiplicity of the eigenvalue ${\mathrm{P}}_{ij}$
of $A_j$ (provided that ${\mathrm{P}}_{ij}\neq {\mathrm{P}}_{kj}$
for $k \neq i$). Clearly, we have $m_0=1$ and $\sum_{i=0}^dm_i=N$
since, $\sum_{i=0}^dE_i=I$. It should be also noticed that, due to
the tracelessness of the adjacency matrices $A_i$, $i\neq0$, one
can obtain the following identity (which we will use it in
obtaining the main result of the paper)
\begin{equation}\label{At}
N\delta_{l0}=trA_l=N.\sum_{j=0}^d
{\mathrm{P}}_{lj}\langle\alpha|E_j|\alpha\rangle=\sum_{j=0}^d
{\mathrm{P}}_{lj}m_j\;\,
\end{equation}
where, we have used the fact that
$\langle\alpha|E_j|\alpha\rangle$ is independent of the choic of
$\alpha\in V$. Finally the underlying network of an association
scheme $\Gamma=(V,R_1)$ is an undirected connected network, where
the set $V$ and $R_1$ consist of its vertices and edges,
respectively. Obviously replacing $R_1$ with one of the other
relations $R_i$, $i\neq 0,1$ will also gives us an underlying
network $\Gamma=(V,R_i)$ (not necessarily a connected network)
with the same set of vertices but a new set of edges $R_i$.

As we will see in subsection $2.3$, in the case of
distance-regular networks, the adjacency matrices $A_j$ are
polynomials of the adjacency matrix $A$, i.e., $A_j=P_j(A)$, where
$P_j$ is a polynomial of degree $j$, then the eigenvalues
${\mathrm{P}}_{ij}$ in (\ref{eign}) are polynomials of the
eigenvalues ${\mathrm{P}}_{i1}\equiv \lambda_i$ (eigenvalues of
the adjacency matrix $A$). This indicates that, in
distance-regular networks (underlying networks of $P$-polynomial
association schemes) the matrix $\mathrm{P}$ is a polynomial
transformation \cite{Puschel}, as
\begin{equation}\label{push}
\mathrm{P}=\left(\begin{array}{ccccc}
    1 & 1 & \ldots & 1 \\
      P_1(\lambda_0) & P_1(\lambda_1) & \ldots & P_1(\lambda_d) \\
      P_2(\lambda_0) & P_2(\lambda_1) & \ldots & P_2(\lambda_d) \\
      \vdots & \vdots & \ldots & \vdots \\
      P_d(\lambda_0) & P_d(\lambda_1) & \ldots & P_d(\lambda_d) \\
    \end{array}\right)
    \end{equation}
or ${\mathrm{P}}_{ij}=P_i(\lambda_j)$. For example, for the
directed cyclic graph $C_n$, where the adjacency matrices are
given by $A_i=S^i$ with $S$ as shift operator of period $n$
($S^n=I$), we have $A_i=A^i$ and so the polynomial transformation
in (\ref{push}) reads as the well known Fourier transformation
(the eigenvalues of the $S$ are the $n$-th roots of unity, i.e.,
$\lambda_j=\omega^i=e^{2\pi j/n}$). For undirected cyclic graph
$C_n$, the adjacency matrices are given by $A_i=S^i+S^{-i}$ and
one can show that $A_i=T_i(A/2)$, where $T_i$'s are Chebyshev
polynomials of the first kind. Then the polynomial transformation
(\ref{push}) is the same as the discrete cosine transformation.
\subsection{Stratifications}
For a given vertex $\alpha\in V$, let
$\Gamma_i(\alpha):=\{\beta\in V: (\alpha, \beta)\in R_i\}$  denote
the set of all vertices having the relation $R_i$ with $\alpha$.
Then, the vertex set $V$ can be written as disjoint union of
$R_i(\alpha)$ for $i=0,1,2,...,d$, i.e.,
 \begin{equation}\label{asso1}
 V=\bigcup_{i=0}^{d}\Gamma_{i}(\alpha).
 \end{equation}
We fix a point $o\in V$ as an origin of the underlying network of
an association scheme, called reference vertex. Then, the relation
(\ref{asso1}) stratifies the network into a disjoint union of
associate classes $\Gamma_{i}(o)$ (called the $i$-th stratum with
respect to $o$). Let $l^2(V)$ denote the Hilbert space of
$C$-valued square-summable functions on $V$. With each associate
class $\Gamma_{i}(o)$ we associate a unit vector in $l^2(V)$
defined by
\begin{equation}\label{unitv}
\ket{\phi_{i}}=\frac{1}{\sqrt{\kappa_i}}\sum_{\alpha\in
\Gamma_{i}(o)}\ket{\alpha},
\end{equation}
where, $\ket{\alpha}$ denotes the eigenket of $\alpha$-th vertex
at the associate class $\Gamma_{i}(o)$ and
$\kappa_i=|\Gamma_{i}(o)|$ is called the $i$-th valency of the
graph ($\kappa_i:=p^0_{ii}=|\{\gamma:(o,\gamma)\in
R_i\}|=|\Gamma_{i}(o)|$).
  The closed subspace of $l^2(V)$ spanned by
$\{\ket{\phi_{i}}\}$ is denoted by $\Lambda(\Gamma)$. Since
$\{\ket{\phi_{i}}\}$ becomes a complete orthonormal basis of
$\Lambda(G)$, we often write
\begin{equation}
\Lambda(\Gamma)=\sum_{i}\oplus \textbf{C}\ket{\phi_{i}}.
\end{equation}
Let $A_i$ be the adjacency matrix of the network $\Gamma=(V,R)$.
Then, from the definition of the $i$-th adjacency matrix $A_i$,
for the reference state $\ket{\phi_0}$ ($\ket{\phi_0}=\ket{o}$,
with $o\in V$ as reference vertex), we have
\begin{equation}\label{Foc1}
A_i\ket{\phi_0}=\sum_{\beta\in \Gamma_{i}(o)}\ket{\beta}.
\end{equation}
Then by using (\ref{unitv}) and (\ref{Foc1}), we have
\begin{equation}\label{Foc2}
A_i\ket{\phi_0}=\sqrt{\kappa_i}\ket{\phi_i}.
\end{equation}
\subsection{Distance-regular graphs}
In this section we consider some set of important graphs called
distance-regular graphs. First we recall the definition of
so-called $P$-polynomial association schemes (which are closely
related to the distance-regular graphs) as follows:\\
\textbf{Definition 2.2} ($P$-polynomial property) The symmetric
association scheme $Y=(X,{\{R_i\}}_{0\leq i\leq d})$ is said to be
$P$-polynomial (with respect to the ordering $R_0,..., R_d$ of the
associate classes) whenever for all $i=0,1,...,d$, there exist
$d_i$, $e_i$, $f_i$; $d_i\neq0\neq f_i$ with:
\begin{equation}\label{P}
A_1A_i=d_iA_{i-1}+e_iA_i+f_iA_{i+1}.
\end{equation}
The condition (\ref{P}) is similar to the well known three term
recursion relations appearing in the theory of orthogonal
polynomials, where $A_1$ is in correspondence with $x$ (see
equation (\ref{trt0}) in subsection $2.4$). Using the recursion
relations (\ref{P}), one can show that $A_i$ is a polynomial in
$A_1$ of degree $i$, i.e., we have $A_i=P_i(A_1)$ for
$i=1,2,...,d$. In particular, $A\equiv A_1$ multiplicatively
generates the Bose-Mesner algebra (for more details see
\cite{bannai}).

An undirected connected graph $\Gamma=(V,R_1)$ is called distance-
regular graph if it is the underlying graph of a $P$-polynomial
association scheme, where the relations are based on distance
function defined as follows: Let the distance between
$\alpha,\beta\in V$ denoted by $\partial(\alpha, \beta)$ is the
length of the shortest walk connecting $\alpha$ and $\beta$
(recall that a finite sequence $\alpha_0, \alpha_1,..., \alpha_n
\in V$ is called a walk of length $n$ if $\alpha_{k-1}\sim
\alpha_k$ for all $k=1, 2,..., n$, where $\alpha_{k-1}\sim
\alpha_k$ means that $\alpha_{k-1}$ is adjacent with
$\alpha_{k}$), then the relations $R_i$ in distance-regular graphs
are defined as: $(\alpha, \beta)\in R_i$ if and only if
$\partial(\alpha, \beta)=i$, for $i=0,1,...,d$, where
$d:=$max$\{\partial(\alpha, \beta): \alpha, \beta\in V \}$ is
called the diameter of the graph. Since $\partial(\alpha, \beta)$
gives the distance between vertices $\alpha$ and $\beta$,
$\partial$ is called the distance function. Clearly, we have
$\partial(\alpha, \alpha)=0$ for all $\alpha\in V$ and
$\partial(\alpha, \beta)=1$ if and only if $\alpha\sim \beta$.
Therefore, distance-regular graphs become metric spaces with the
distance function $\partial$.

One should notice that, the condition (\ref{P}) implies that for
distance-regular graphs, we have the following relation
\begin{equation}\label{reldis}
\Gamma_1(\beta)\subseteq \Gamma_{i-1}(\alpha)\cup
\Gamma_i(\alpha)\cup \Gamma_{i+1}(\alpha),\;\;\ \forall \;\
\beta\in \Gamma_i(\alpha).
\end{equation}
We also note that, in distance-regular graphs, the stratification
is reference vertex independent, namely one can choose every
vertex as a reference one, while the stratification of more
general graphs may be reference dependent.

The relation (\ref{P}) implies that in a distance-regular graph,
$p_{j1}^i=0 $ (for $i\neq 0$, $j$ dose not belong to $\{i-1, i,
i+1 \}$), i.e., the non-zero intersection numbers of the graph are
given by
\begin{equation}\label{abc}
 a_i=p_{i1}^i, \;\;\;\  b_i=p_{i+1,1}^i, \;\;\;\
 c_i=p_{i-1,1}^i\;\ ,
\end{equation}
respectively (see Figure $3$). The intersection numbers
(\ref{abc}) and the valencies $\kappa_i$ satisfy the following
obvious conditions
$$a_i+b_i+c_i=\kappa,\;\;\ \kappa_{i-1}b_{i-1}=\kappa_ic_i ,\;\;\
i=1,...,d,$$
\begin{equation}\label{intersec}
\kappa_0=c_1=1,\;\;\;\ b_0=\kappa_1=\kappa, \;\;\;\ (c_0=b_d=0).
\end{equation}
Thus all parameters of the graph can be obtained from the
intersection array $\{b_0,...,b_{d-1};c_1,...,c_d\}$.

By using the equations (\ref{ss}) and (\ref{intersec}), for
adjacency matrices of distance-regular graph $\Gamma$, we obtain
$$A_1A_i=b_{i-1}A_{i-1}+a_iA_i+c_{i+1}A_{i+1},\;\
i=1,2,...,d-1,$$
\begin{equation}\label{P0}
A_1A_d=b_{d-1}A_{d-1}+(\kappa-c_d)A_d.
\end{equation}
The recursion relations (\ref{P0}), imply that
\begin{equation}\label{P1}
A_i=P_i(A),\;\ i=0,1,...,d.
\end{equation}
By acting two sides of (\ref{P0}) on $\ket{\phi_0}$ and using
(\ref{Foc2}), we obtain
\begin{equation}\label{P2}
\sqrt{\kappa_i}A\ket{\phi_i}=\sqrt{\kappa_{i-1}}b_{i-1}\ket{\phi_{i-1}}+
\sqrt{\kappa_i}a_i\ket{\phi_i}+\sqrt{\kappa_{i+1}}c_{i+1}\ket{\phi_{i+1}},\;\
i=0,1,...,d.
\end{equation}
Then, by dividing the sides of (\ref{P2}) by $\sqrt{\kappa_i}$ and
using (\ref{intersec}), one can easily obtain the following three
term recursion relations for the unit vectors $\ket{\phi_i}$,
$i=0,1,...,d$
\begin{equation}\label{trt}
A\ket{\phi_i}=\beta_{i+1}\ket{\phi_{i+1}}+\alpha_i\ket{\phi_i}+\beta_{i}\ket{\phi_{i-1}},
\end{equation}
where, the coefficients $\alpha_i$ and $\beta_i$ are defined
as
\begin{equation}\label{omegal}
\alpha_0=0,\;\;\ \alpha_k\equiv a_k=\kappa-b_{k}-c_{k},\;\;\;\;\
\omega_k\equiv\beta^2_k=b_{k-1}c_{k},\;\;\ k=1,...,d.
\end{equation}
That is, in the basis of the unit vectors
$\{\ket{\phi_i},i=0,1,...,d\}$, the adjacency matrix $A$ is
projected to the following symmetric tridiagonal form:
\begin{equation}\label{trid}
A=\left(
\begin{array}{cccccc}
 \alpha_0 & \beta_1 & 0 & ... &...&0 \\
      \beta_1 & \alpha_1 & \beta_2 & 0 &...&0 \\
      0 & \beta_2 & \alpha_3 & \beta_3 & \ddots&\vdots \\
     \vdots & \ddots &\ddots& \ddots &\ddots &0\\
     0 & \ldots  &0 &\beta_{d-1} & \alpha_{d-1} &\beta_{d}\\
          0&... & 0 &0 & \beta_{d} & \alpha_{d}\\
\end{array}
\right).
\end{equation}
In Ref. \cite{js2}, it has be shown that, the coefficients
$\alpha_i$ and $\beta_i$ can be also obtained easily by using the
Lanczos iteration algorithm.
\subsection{Stieltjes function associated with the network }
In this subsection we recall the definition of the Stieltjes
function associated with a distance-regular network which is
related to the spectral distribution corresponding to the network.
To do so, first we recall some facts about the spectral
distribution associated with the adjacency matrix of the network.
In fact, the spectral analysis of operators is an important issue
in quantum mechanics, operator theory and mathematical physics
\cite{simon, Hislop}. Since the advent of random matrix theory
(RMT), there has been considerable interest in the statistical
analysis of spectra \cite{rmt,rmt1,rmt3}. RMT can be viewed as a
generalization of the classical probability calculus, where the
concept of probability density distribution for a one-dimensional
random variable is generalized onto an averaged spectral
distribution of the ensemble of large, non-commuting random
matrices. Such a structure exhibits several phenomena known in
classical probability theory, including central limit theorems
\cite{centl}. Also, the two-point resistance has a probabilistic
interpretation based on classical random walker walking on the
network. Indeed, the connection between random walks and electric
networks has been recognized for some time (see e.g.
\cite{Kakutani, Kemeny, Kelly} ), where one can establish a
connection between the electrical concepts of current and voltage
and corresponding descriptive quantities of random walks regarded
as finite state Markov chains (for more details see \cite{Doyle}).
Also, by adapting the random-walk dynamics and mean-field theory
it has been studied that \cite{Bosiljka}, how the growth of a
conducting network, such as electrical or electronic circuits,
interferes with the current flow through the underlying evolving
graphs. In Ref.\cite{js2} it has been shown that, there is also
connection between the mathematical techniques such as Hilbert
space of the stratification and spectral techniques (which have
been employed in \cite{js,js1,rootlatt,laplac} for investigating
continuous time quantum walk on graphs), and electrical concept of
resistance between two arbitrary nodes of regular networks and so
the same techniques can be used for calculating the resistance.
Note that, although we take the spectral approach to define the
Stieltjes function in terms of orthogonal polynomials (which are
orthogonal with respect to the spectral distribution $\mu$
associated with the network) with three term recursion relations,
in practice as it will be seen in the section $3$, we will
calculate two-point resistances without any need to evaluate the
spectral distribution $\mu$.

It is well known that, for any pair $(A,\ket{\phi_0})$ of a matrix
$A$ and a vector $\ket{\phi_0}$, it can be assigned a measure
$\mu$ as follows
\begin{equation}\label{sp1}
\mu(x)=\braket{ \phi_0}{E(x)|\phi_0},
\end{equation}
 where
$E(x)=\sum_i|u_i\rangle\langle u_i|$ is the operator of projection
onto the eigenspace of $A$ corresponding to eigenvalue $x$, i.e.,
\begin{equation}
A=\int x E(x)dx.
\end{equation}
It is easy to see that, for any polynomial $P(A)$ we have
\begin{equation}\label{sp2}
P(A)=\int P(x)E(x)dx,
\end{equation}
where for discrete spectrum the above integrals are replaced by
summation. Therefore, using the relations (\ref{sp1}) and
(\ref{sp2}), the expectation value of powers of adjacency matrix
$A$ over starting site $\ket{\phi_0}$ can be written as
\begin{equation}\label{v2}
\braket{\phi_{0}}{A^m|\phi_0}=\int_{R}x^m\mu(dx), \;\;\;\;\
m=0,1,2,....
\end{equation}
The existence of a spectral distribution satisfying (\ref{v2}) is
a consequence of Hamburger's theorem, see e.g., Shohat and
Tamarkin [\cite{st}, Theorem 1.2].

Obviously relation (\ref{v2}) implies an isomorphism from the
Hilbert space of the stratification onto the closed linear span of
the orthogonal polynomials with respect to the measure $\mu$. More
clearly, the orthonormality of the unit vectors $\ket{\phi_i}$
implies that
\begin{equation}\label{ortpo}
\delta_{ij}=\langle\phi_i|\phi_j\rangle=\frac{1}{\sqrt{\kappa_i\kappa_j}}\braket{\phi_{0}}{A_iA_j|\phi_0}=\int_{R}P'_i(x)P'_j(x)\mu(dx),
\end{equation}
where, we have used the equations (\ref{Foc2}) and (\ref{P1}) to
write
\begin{equation}\label{xx}
\ket{\phi_i}=\frac{1}{\sqrt{\kappa_i}}A_i\ket{\phi_0}=\frac{1}{\sqrt{\kappa_i}}P_i(A)\ket{\phi_0}\equiv
P'_i(A)\ket{\phi_0},
\end{equation}
with $P'_i(A):=\frac{1}{\sqrt{\kappa_i}}P_i(A)$. Now, by
substituting (\ref{xx}) in (\ref{trt}), we get three term
recursion relations between polynomials $P'_j(A)$, which leads to
the following  three term recursion relations between polynomials
$P'_j(x)$
\begin{equation}\label{trt0}
xP'_{k}(x)=\beta_{k+1}P'_{k+1}(x)+\alpha_kP'_{k}(x)+\beta_kP'_{k-1}(x)
\end{equation}
for $k=0,...,d-1$, with $P'_0(x)=1$. Multiplying two sides of
(\ref{trt0}) by $\beta_1...\beta_k$ we obtain
\begin{equation}
\beta_1...\beta_kxP'_{k}(x)=\beta_1...\beta_{k+1}P'_{k+1}(x)+\alpha_k\beta_1...\beta_kP'_{k}(x)+\beta_k^2.\beta_1...\beta_{k-1}P'_{k-1}(x).
\end{equation}
By rescaling $P'_k$ as $Q_k=\beta_1...\beta_kP'_k$, the spectral
distribution $\mu$ under question is characterized by the property
of orthonormal polynomials $\{Q_k\}$ defined recurrently by
$$ Q_0(x)=1, \;\;\;\;\;\
Q_1(x)=x,$$
\begin{equation}\label{op}
xQ_k(x)=Q_{k+1}(x)+\alpha_{k}Q_k(x)+\beta_k^2Q_{k-1}(x),\;\;\
k\geq 1.
\end{equation}
It is a well known result from the theory of orthogonal
polynomials that, the polynomials defined by (\ref{op}) can be
also evaluated via the following determinant
\begin{equation}\label{op'}
Q_k(x)=\left|\begin{array}{ccccccccc}
             x-\alpha_0 & 1 & 0 & 0 & \ldots & 0 & 0\;\ & 0\;\ & 0 \\
             \omega_1\;\ & x-\alpha_1 & 1 & 0 & \ldots & 0 & 0\;\ & 0\;\ & 0 \\
             0 \;\ & \omega_2\;\ & x-\alpha_2 & 1 & \ldots & 0 & 0\;\ & 0\;\ & 0 \\
             \vdots \;\ & \vdots \;\ & \vdots & \vdots & \ldots & \vdots & \vdots\;\ & \vdots\;\ & \vdots \\
             0 \;\ & 0\;\ & 0 & 0 & \ldots & 0 & \omega_{k-2}\;\ & x-\alpha_{k-2} & 1 \\
             0 \;\ & 0\;\ & 0 & 0 & \ldots & 0 & 0\;\ & \omega_{k-1}\;\ & x-\alpha_{k-1} \\
           \end{array}\right|
\end{equation}

If such a spectral distribution is unique, the spectral
distribution $\mu$ is determined by the identity
\begin{equation}\label{sti}
G_{\mu}(x)=\int_{R}\frac{\mu(dy)}{x-y}=\frac{1}{x-\alpha_0-\frac{\omega_1}{x-\alpha_1-\frac{\omega_2}
{x-\alpha_2-\frac{\omega_3}{x-\alpha_3-\cdots}}}}=\frac{Q_{d-1}^{(1)}(x)}{Q_{d}(x)}=\sum_{l=0}^{d-1}
\frac{b_l}{x-x_l},
\end{equation}
where, $x_l$ are the roots of the polynomial $Q_{d}(x)$. The
function $G_{\mu}(x)$ is called the Stieltjes/Hilbert transform of
spectral distribution $\mu$ or Stieltjes function and polynomials
$\{Q_{k}^{(1)}\}$ are defined recurrently as
$$Q_{0}^{(1)}(x)=1, \;\;\;\;\;\
    Q_{1}^{(1)}(x)=x-\alpha_1,$$
\begin{equation}\label{oq}
xQ_{k}^{(1)}(x)=Q_{k+1}^{(1)}(x)+\alpha_{k+1}Q_{k}^{(1)}(x)+\beta_{k+1}^2Q_{k-1}^{(1)}(x),\;\;\
k\geq 1,
\end{equation}
respectively. The coefficients $b_l$ appearing in (\ref{sti}) are
calculated as
\begin{equation}\label{Gauss}
b_l=\lim_{x\rightarrow x_l}(x-x_l)G_{\mu}(x).
\end{equation}

Now if $G_{\mu}(x)$ is known, then the spectral distribution $\mu$
can be recovered from $G_{\mu}(x)$ by means of the Stieltjes
inversion formula:
\begin{equation}\label{oq1}
\mu(y)- \mu(x)=-\frac{ 1}{\pi} lim_{ v\rightarrow 0^+}\int_{x}^{y}
Im{G_{\mu}(u+ iv)}du.
\end{equation}
Substituting the right hand side of (\ref{sti}) in (\ref{oq1}),
the spectral distribution can be determined in terms of $x_l$, $l
= 1, 2, ...,$ the roots of the polynomial $Q_d(x)$, and Guass
quadrature constants $b_l$, $l = 1, 2, ...$ as
\begin{equation}\label{oq2}
\mu(x)=\sum_{l}b_l\delta(x - x_l)
\end{equation}
(for more details see Refs.\cite{st, obh,tsc,obah}). Note that, by
using (\ref{sp1}) and (\ref{oq2}), we have
\begin{equation}\label{oq3}
b_i=\mu(x_i)=\langle\alpha|E(x_i)|\alpha\rangle=\frac{1}{N}tr
E_i=\frac{m_i}{N}.
\end{equation}
\subsubsection{The Christoffel-Darboux Identity}
In the next section we will treat with the calculation of the
effective resistances in distance-regular networks, where we will
use one of the most important theorems of the orthogonal
polynomials known as Christoffel-Darboux Identity. This identity
is expressed as follows\\
\textbf{Theorem (Christoffel-Darboux Identity)} Let $\{Q_n(x)\}$
satisfy (\ref{op}). Then
\begin{equation}\label{CH-D}
\sum_{k=1}^n\frac{Q_k(x)Q_k(u)}{\omega_1\omega_2...\omega_k}=(\omega_1\omega_2...\omega_n)^{-1}\frac{Q_{n+1}(x)Q_n(u)-Q_n(x)Q_{n+1}(u)}{x-u}.
\end{equation}
For the proof, the reader is refered to \cite{tsc}.
\section{Two-point resistances in regular resistor networks}
A classic problem in electric circuit theory studied by numerous
authors over many years, is the computation of the resistance
between two nodes in a resistor network (see, e.g.,
\cite{Cserti}). The results obtained in this section show that,
there is a close connection between the techniques introduced in
section $2$ such as Hilbert space of the stratification and
electrical concept of resistance between two arbitrary nodes of
regular networks and these techniques can be employed for
calculating two-point resistances.

For a given regular graph $\Gamma$ with $v$ vertices and adjacency
matrix $A$, let $r_{ij}=r_{ji}$ be the resistance of the resistor
connecting vertices $i$ and $j$. Hence, the conductance is
$c_{ij}=r^{-1}_{ij}=c_{ji}$ so that $c_{ij}=0$ if there is no
resistor connecting $i$ and $j$. Denote the electric potential at
the $i$-th vertex by $V_i$ and the net current flowing into the
network at the $i$-th vertex by $I_i$ (which is zero if the $i$-th
vertex is not connected to the external world). Since there exist
no sinks or sources of current including the external world, we
have the constraint $\sum_{i=1}^vI_i=0$. The Kirchhoff law states
\begin{equation}\label{resistor}
\sum_{j=1,j\neq i}^vc_{ij}(V_i-V_j)=I_i,\;\;\  i=1,2,...,v.
\end{equation}
Explicitly, Eq.(\ref{resistor}) reads
\begin{equation}\label{resistor1}
L\vec{V}=\vec{I},
\end{equation}
where, $\vec{V}$ and $\vec{I}$ are $v$-vectors whose components
are $V_i$ and $I_i$, respectively and
\begin{equation}\label{laplas}
L=\sum_{i}c_i|i\rangle\langle i|-\sum_{i,j}c_{ij}|i\rangle\langle
j|
\end{equation}
is the Laplacian of the graph $\Gamma$ with
\begin{equation}
c_i\equiv \sum_{j=1,j\neq i}^vc_{ij},
\end{equation}
for each vertex $\alpha$. Hereafter, we will assume that all nonzero
resistances are equal to $1$, then the off-diagonal elements of $-L$
are precisely those of $A$, i.e.,
\begin{equation}\label{laplas1}
L=\kappa I-A,
\end{equation}
with $\kappa=deg(\alpha)$, for each vertex $\alpha$. It should be
noticed that, $L$ has eigenvector $(1,1,...,1)^t$ with eigenvalue
$0$. Therefore, $L$ is not invertible and so we define the
psudo-inverse of $L$ as
\begin{equation}\label{inv.laplas}
L^{-1}=\sum_{i,\lambda_i\neq0} {\lambda}^{-1}_iE_i,
\end{equation}
where, $E_i$ is the operator of projection onto the eigenspace of
$L^{-1}$ corresponding to eigenvalue $\lambda_i$. It has been
shown that, the two-point resistances $R_{\alpha\beta}$ are given
by
\begin{equation}\label{eq.res.}
R_{\alpha\beta}=\langle \alpha|L^{-1}|\alpha\rangle+\langle
\beta|L^{-1}|\beta\rangle-\langle
\alpha|L^{-1}|\beta\rangle-\langle \beta|L^{-1}|\alpha\rangle.
\end{equation}
This formula may be formally derived using Kirchoff 's laws, and
seems to have been long known in the electrical engineering
literature, with it appearing in several texts, such as
Ref.\cite{12}. For distance-regular graphs as resistor networks,
the diagonal entries of $L^{-1}$ are independent of the vertex,
i.e., $L^{-1}_{\alpha\alpha}=L^{-1}_{\beta\beta}$ for all
$\alpha,\beta\in V$. Therefore, from the relation (\ref{eq.res.})
and the fact that $L^{-1}$ is a real matrix, we can obtain the
two-point resistance between two arbitrary nodes $\alpha$ and
$\beta$ as follows
\begin{equation}\label{eq.res.dist.}
R_{\alpha\beta}=2(L^{-1}_{\alpha\alpha}-L^{-1}_{\alpha\beta}).
\end{equation}
It should be noticed that, in distance-regular graphs (more
generally, for underlying networks of association scheme) due to
the stratification of the network, all of the nodes belonging to
the same stratum with respect to the reference node, i.e.,
$\alpha$, possess the same two-point resistance with respect to
the $\alpha$. More clearly, for all $\beta\in \Gamma_m(\alpha)$ we
have
\begin{equation}\label{eq.res.dist.1}
L^{-1}_{\alpha\beta^{(m)}}=\langle\alpha|L^{-1}|\beta\rangle=\frac{1}{\sqrt{\kappa_m}}\langle\alpha|L^{-1}|\phi_m\rangle=\frac{1}{\kappa_m}\langle\alpha|A_mL^{-1}|\alpha\rangle.
\end{equation}
Then, by using (\ref{eq.res.dist.}), we obtain
\begin{equation}\label{eq.res.dist.2}
R_{\alpha\beta^{(m)}}=\frac{2}{\kappa_m}\{\kappa_mL^{-1}_{\alpha\alpha}-(A_mL^{-1})_{\alpha\alpha}\}=\frac{2}{\kappa_m}\langle\alpha|\frac{\kappa_m1-P_m(A)}{\kappa1-A}|\alpha\rangle,\;\
\forall \;\ \beta\in \Gamma_m(\alpha).
\end{equation}
where, the upper index $m$ in $L^{-1}_{\alpha\beta^{(m)}}$ and
$R_{\alpha\beta^{(m)}}$ indicate that $\beta$ belongs to the
$m$-th stratum with respect to $\alpha$.

As, Eq.(\ref{eq.res.dist.}) implies, in order to evaluate
two-point resistances $R_{\alpha\beta^{(m)}}$, we need to
calculate the matrix entries $L^{-1}_{\alpha\alpha}$ and
$L^{-1}_{\alpha\beta}$. To this end, one needs to know the
spectrum of the pseudo inverse $L^{-1}$ (see
Eq.(\ref{inv.laplas})) which is a task with high complexity for
networks with large size, even with computer. In the following we
give an algebraic method such that the two-point resistances are
calculated recursively without any knowledge about the spectrum of
the pseudo inverse of Laplacian of the network.
\subsection{Random walks and electrical networks}
The computation of effective resistances is relevant to a wide
range of problems ranging from random walks (see Ref.
\cite{Doyle}). Random walks on graphs are the bases of a number of
classical algorithms. Examples include 2-SAT (satisfiability for
certain types of Boolean formulas), graph connectivity, and
finding satisfying assignments for Boolean formulas. It is this
success of random walks that motivated the study of their quantum
analogs in order to explore whether they might extend the set of
quantum algorithms. In Refs. \cite{js,js1,rootlatt,js2,laplac},
the same techniques introduced in this paper in order to evaluate
the effective resistances, such as the algebraic structure of
distance-regular graphs (Bose-Mesner algebra), stratification and
spectral analysis methods have been used for investigation of the
continous time quantum walks on the regular networks. In order to
show this connection more clearly, in the following we discuss the
link between resistance distance and two important quantities
(average first passage time and average commute time) defined in
random walks on graphs.

Let $\Gamma$ be a complete undirected graph with $N$ vertices
numbered $1,2,. . .,N$, in which each edge $(\alpha,\beta)$ is
assigned its distance $\partial(\alpha,\beta)\equiv d_{\alpha\beta}
> 0$. One can study the harmonic random walk in $\Gamma$ with escape
probability (the probability that a walk starting at $\alpha$
reaches $\beta$ before it returns to $\alpha$)
\begin{equation}\label{pesc}
p_{esc}(\alpha,\beta)=\frac{1/d_{\alpha\beta}}{\sum_{\gamma\neq\alpha}1/d_{\alpha\gamma}}.
\end{equation}

In the literature, the harmonic random walk was often studied using
techniques from electrical network theory.  Denote by
$H_{\alpha\beta}$, the hitting cost of the harmonic walk from
$\alpha$ to $\beta$, defined as the expected cost (total distance)
to reach $\beta$ for the first time when started from $\alpha$. By
elementary probability, the costs of the random walks ending at
vertex $N$ satisfy the following system of equations
\begin{equation}\label{eqx}
H_{\alpha\beta}=\sum_{\beta\neq \alpha}
p_{\alpha\beta}(d_{\alpha\beta} + H_{\beta,N}),\;\ \mathrm{for}\;\
\alpha\neq N,
\end{equation}
and $H_{N,N}=0$. As noted in \cite{Doyle}, the hitting costs have
an interpretation in terms of electrical networks. We can think of
$\Gamma$ as an electrical network in which each edge
$(\alpha,\beta)$ has resistance $d_{\alpha\beta}$ (so, the
Eq.(\ref{pesc}) is a probabilistic interpretation of the effective
conductance). If we inject current of value $N-1$ into each node
and draw current of value $N(N -1)$ from node $N$, then the
voltages relative to node $N$ established at the nodes satisfy the
same equation as (\ref{eqx}). Therefore the voltage at $\alpha$ is
equal to the hitting cost $H_{\alpha\beta}$. Moreover, for a
random walker on a network, one can assign a quantity known as
average commute time (CT) denoted by $n(\alpha, \beta)$ which is
defined as the average number of steps the random walker, starting
in state $\alpha\neq\beta$, will take before entering a given
state $\beta$ for the first time, and go back to $\alpha$.
Clearly, the average commute time is symmetric and is equal to
$n(\alpha, \beta)=m(\beta|\alpha)+m(\alpha|\beta)$, where
$m(\beta|\alpha)$ is the average first-passage time defined as the
average number of steps the random walker, starting in state
$\alpha$, will take to enter state $\beta$ for the first time.
Note that for random walks on distance-regular graphs, we have
$m(\beta|\alpha)=m(\alpha|\beta)$ and so, $n(\alpha,
\beta)=2m(\beta|\alpha)$. By viewing the graph as an electrical
network, the average commute time has an electrical equivalent
\begin{equation}\label{CT}
n(\alpha, \beta)=N.\kappa R_{\alpha\beta},
\end{equation}
where, $N.\kappa$ is the volume of the graph (the volume of a graph
is defined as $\sum_{\alpha\in V}d_{\alpha}$, with $d_{\alpha}$ as
the degree of the vertex $\alpha$). The Eq.(\ref{CT}) indicates
that, the average commute time and effective resistance basically
measure the same quantity. This quantity can also be called
resistance distance (it has be shown that $n(\alpha, \beta)$ is a
distance measure). Further connections between random walks and
electrical networks were explored by Chandra et al.\cite{3new}. In
the following we introduce a method for recursive calculation of
resistance distances on distance-regular resistor networks based on
spectral techniques specially by employing the Christoffel-Darboux
identity.
\section{Recursive calculation of resistance distance based on spectral analysis methods and Christoffel-Darboux identity}
In this section, we show that the resistance distance on
distance-regular resistor networks increases with the number of
the strata, i.e., $R_{\alpha\beta^{(m+1)}}$ is strictly larger
than $R_{\alpha\beta^{(m)}}$ for $m=1,2,...,d-1$. In previous work
\cite{jss1}, explicit formulas for the two-point resistances up to
the third stratum, i.e., $R_{\alpha\beta^{(m)}}$ for $m=1,2,3$,
have been given in terms of the intersection array of the network,
where the authors have been employed the properties of the
Stieltjes function associated with the network. Here in this work,
we use the spectral techniques and Christoffel-Darboux identity
and give a recursive formula for calculating the two-point
resistances $R_{\alpha\beta^{(m)}}$ for $m=1,2,...,d$. In addition
to the preference that this formula enables us to calculate
recursively all of the resistance distances, it also indicates
that the resistance distance strictly increases with the number of
the stratum.

Let $\alpha$ and $\beta$ be two arbitrary nodes of the network
such that $\beta$ belongs to the $m$-th stratum with respect to
$\alpha$, i.e., $\beta\in \Gamma_m(\alpha)$ (we choose one of the
nodes, here $\alpha$, as reference node). Now, for calculating the
matrix entries $L^{-1}_{\alpha\alpha}$ and $L^{-1}_{\beta\alpha}$
in (\ref{eq.res.}), we use the spectral techniques to obtain
\begin{equation}\label{l1}
L^{-1}_{\alpha\alpha}=\langle \alpha|\frac{1}{\kappa
I-A}|\alpha\rangle=\sum_{i=1}^{d}\langle \alpha|\frac{E_i}{\kappa
-\lambda_i}|\alpha\rangle=\frac{1}{N}\sum_{i=1}^{d}\frac{m_i}{\kappa
-\lambda_i},
\end{equation}
and
$$L^{-1}_{\beta\alpha}=\langle \beta|\frac{1}{\kappa
I-A}|\alpha\rangle=\frac{1}{\sqrt{\kappa_m}}\langle
\phi_m|\frac{1}{\kappa
I-A}|\alpha\rangle=\frac{1}{\sqrt{\kappa_m}}\langle
\alpha|\frac{P'_m(A)}{\kappa I-A}|\alpha\rangle=$$
\begin{equation}\label{l2}
\frac{1}{\kappa_m}\sum_{i=1}^d\langle
\alpha|\frac{P_m(\lambda_i)E_i}{\kappa
-\lambda_i}|\alpha\rangle=\frac{1}{N\kappa_m}\sum_{i=1}^{d}\frac{m_iP_m(\lambda_i)}{\kappa-\lambda_i},
\end{equation}
where, we have considered $\lambda_0=\kappa$ ($\kappa$ is the
eigenvalue corresponding to the idempotent $E_0$). Then, by using
Eq.(\ref{eq.res.dist.}) and the fact that $L_{\alpha\alpha}^{-1}$
is independent of $m$ (the number of stratum), we have
$$\hspace{-2cm}R_{\alpha\beta^{(m+1)}}-R_{\alpha\beta^{(m)}}=2(L^{-1}_{\alpha\beta^{(m)}}-L^{-1}_{\alpha\beta^{(m+1)}})=\frac{2}{N}\sum_{j=1}^d\frac{m_j}{\kappa-\lambda_j}(\frac{P_m(\lambda_j)}{\kappa_m}-\frac{P_{m+1}(\lambda_j)}{\kappa_{m+1}})=$$
\begin{equation}\label{res1}
\hspace{-1.5cm}\small{{\frac{2}{N\kappa_m\kappa_{m+1}}\sum_{j=1}^dm_j\frac{\kappa_{m+1}P_{m}(\lambda_j)-\kappa_mP_{m+1}(\lambda_j)}{\kappa-\lambda_j}=\frac{2\sqrt{\kappa_m\kappa_{m+1}}}{N\kappa_m\kappa_{m+1}\beta^2_1...\beta^2_m\beta_{m+1}}\sum_{j=1}^dm_j\frac{Q_{m+1}(\kappa)Q_{m}(\lambda_j)-Q_{m}(\kappa)Q_{m+1}(\lambda_j)}{\kappa-\lambda_j}}}.
\end{equation}
Now, from the fact that $\{Q_k(x)\}$ satisfy the three-term
recursion relations (\ref{op}), we can use the Christoffel-Darboux
identity to write the right hand side of (\ref{res1}) as follows
$$\hspace{-2cm}r.h.s=\frac{2}{N\sqrt{\kappa_m\kappa_{m+1}\omega_{m+1}}.\omega_1...\omega_m}\sum_{j=1}^dm_j(\omega_1...\omega_m\sum_{l=0}^m\frac{Q_l(\kappa)Q_l(\lambda_j)}{\omega_1...\omega_l})=$$
$$\frac{2}{N\sqrt{\kappa_m\kappa_{m+1}\omega_{m+1}}}\sum_{l=0}^m\frac{Q_l(\kappa)}{\omega_1...\omega_l}.\sum_{j=1}^dm_jQ_l(\lambda_j)=
\frac{2}{N\sqrt{\kappa_m\kappa_{m+1}\omega_{m+1}}}\sum_{l=0}^m\frac{Q_l(\kappa)\beta_1...\beta_l}{\omega_1...\omega_l\sqrt{\kappa_l}}.\sum_{j=1}^dm_jP_l(\lambda_j)=$$
\begin{equation}\label{res1'}
\frac{2}{N\sqrt{\kappa_m\kappa_{m+1}\omega_{m+1}}}\sum_{l=0}^m\sum_{j=1}^dm_jP_l(\lambda_j)=\frac{2}{N\sqrt{\kappa_m\kappa_{m+1}\omega_{m+1}}}\sum_{l=0}^m(N\delta_{l0}-\kappa_l),
\end{equation}
where, we have used the fact that
$$Q_l(\kappa)=\frac{\beta_1\beta_2...\beta_l}{\sqrt{\kappa_l}}P_l(\kappa)=\frac{\beta_1\beta_2...\beta_l}{\sqrt{\kappa_l}}\kappa_l,\;\ \mathrm{and}$$
and have done the following simplification by using
(\ref{intersec})
$$\kappa_m\kappa_{m+1}\omega_{m+1}=\kappa_m\kappa_{m+1}b_mc_{m+1}=(b_m\kappa_m)^2.$$
In the last equality of (\ref{res1'}), we have used the
distance-regularity of the network to substitute
$P_l(\lambda_j)={\mathrm{P}}_{ij}$ and then use the Eq.(\ref{At})
(recall that $m_0=1$ and $P_l(\lambda_0)=P_l(\kappa)=\kappa_l$).
After these simplifications, we obtain the main result of the
paper as follows
\begin{equation}\label{res3}
R_{\alpha\beta^{(m+1)}}-R_{\alpha\beta^{(m)}}=\frac{2}{N\kappa_mb_m}(N-\sum_{l=0}^m\kappa_l)>0,\;\;\
m=1,2,...,d-1.
\end{equation}

It should noticed that, by using the formula (\ref{res3}) we can
evaluate the effective resistance between any two nodes
recursively, if we know the two-point resistance
$R_{\alpha\beta^{(1)}}$. In order to calculate
$R_{\alpha\beta^{(1)}}$, one can use the spectral techniues
introduced in subsection $2.4$ (equations (\ref{v2}),(\ref{oq2})
and (\ref{oq3})) to write
\begin{equation}\label{l1}
L^{-1}_{\alpha\alpha}=\langle \alpha|\frac{1}{\kappa
I-A}|\alpha\rangle=\int_{R-\{\kappa\}}\frac{d\mu(x)}{\kappa-x}=\frac{1}{N}\sum_{i,i\neq0}^{d-1}\frac{m_i}{\kappa-x_i}
\end{equation}
and
\begin{equation}\label{l2}
L^{-1}_{\alpha\beta^{(1)}}=\frac{1}{\kappa}\langle
\alpha|\frac{A}{\kappa I-A}|\alpha\rangle
=\frac{1}{\kappa}\int_{R-\{\kappa\}}\frac{d\mu(x)}{\kappa-x}x=\frac{1}{N\kappa}\sum_{i,i\neq0}\frac{m_ix_i}{\kappa-x_i},
\end{equation}
Then, by using (\ref{eq.res.dist.}), one can obtain
\begin{equation}\label{l3}
R_{\alpha\beta^{(1)}}=\frac{2}{N\kappa}\sum_{i,i\neq0}\frac{m_i(\kappa-x_i)}{\kappa-x_i}=\frac{2}{N\kappa}\sum_{i,i\neq0}m_i=\frac{2(N-1)}{N\kappa}.
\end{equation}

Note that, as the main result of the paper, the result
(\ref{res3}) shows that the resistance distance (and consequently
the other quantities related to the resistance distance such as
the average first passage time and Euclidean commute time
associated with a random walk) on distance-regular networks is
strictly increasing function of the shortest path distance defined
in subsection $2.3$, i.e., the nodes belonging to the farthest
stratum with respect to $\alpha$, possess the smallest effective
resistance with $\alpha$. Apart from this fact, the formula
(\ref{res3}) together with (\ref{l3}) gives an algebraic method
for calculation of the resistance distances on distance-regular
networks, where one needs only to know the intersection array of
the networks without any knowledge about the spectrum of the
pseudo inverse of Laplacian of the networks.
\section{Examples}
In this section, we calculate the effective resistances on the
examples of Cycle network, $d$-cube network and Johnson network
recursively by using the formula (\ref{res3}), where for the first
example (cycle network), the general formula for the effective
resistances $R_{\alpha\beta^{(i)}}$, $i=1,2,...,d$ is deduced,
whereas for the two latter ones the effective resistances are
calculated up to the third stratum. The effective resistances
$R_{\alpha\beta^{(i)}}$, $i=1,2,...,d$ on some other important
finite distance-regular networks is given in the appendix.
\subsection{Cycle network $C_{N}$}
The graph $C_N$ for $N=2m$ or $N=2m+1$ consists of $m+1$ strata. The
intersection arrays for even and odd number of vertices is given by
\begin{equation}\label{intcye.}
\{b_0,...,b_{m-1};c_1,...,c_m\}=\{2,1,...,1,1;1,...,1,2\}
\end{equation}
\begin{equation}\label{intcyo.}
\mbox{and}\;\
\{b_0,...,b_{m-1};c_1,...,c_m\}=\{2,1,...,1;1,...,1,1\},
\end{equation}
respectively. We consider the even case $N=2m$, the odd case can be
considered similarly. For this case, we have $\kappa_0=\kappa_m=1;$
$\kappa_l=2$, for $l=1,...,m-1$. Then by using (\ref{l3}) and
(\ref{res3}), we obtain the effective resistances recursively as
follows
$$R_{\alpha\beta^{(1)}}=\frac{2m-1}{2m},$$
$$R_{\alpha\beta^{(2)}}=R_{\alpha\beta^{(1)}}+\frac{2m-3}{2m}=\frac{2(2m-2)}{2m},$$
\begin{equation}\label{eqcycle}
R_{\alpha\beta^{(3)}}=R_{\alpha\beta^{(2)}}+\frac{2m-5}{2m}=\frac{3(2m-3)}{2m},...
\end{equation}
From (\ref{eqcycle}), one can easily deduce the following result
\begin{equation}\label{eqcycle1}
R_{\alpha\beta^{(l)}}=\frac{l(2m-l)}{2m},\;\;\ l=1,2,...,m.
\end{equation}

The above formula indicates that in the limit of the large number of
vertices, i.e., in the limit $m\rightarrow \infty$, where the cycle
network tends to the infinite line network, the effective
resistances are given by
\begin{equation}\label{eqinfline}
R_{\alpha\beta^{(l)}}=l,\;\;\ l=1,2,... \quad.
\end{equation}
\subsection{$d$-cube}
The $d$-cube, i.e. the hypercube of dimension $d$, also called
Hamming cube, is a network with $2^d$ nodes, each of which can be
labeled by an $d$-bit binary string. Two nodes on the hypercube
described by bitstrings $\vec{x}$ and $\vec{y}$ are are connected by
an edge if $|\vec{x}- \vec{y}|=1$, where $|\vec{x}|$ is the Hamming
weight of $\vec{x}$. In other words, if $\vec{x}$ and $\vec{y}$
differ by only a single bit flip, then the two corresponding nodes
on the graph are connected. Thus, each of the $2^d$ nodes on the
$d$-cube has degree $d$. For the $d$-cube we have $d+1$ strata with
\begin{equation}\label{cubevalenc}
\kappa_i=\frac{d!}{i!(d-i)!}\;\ , \;\ 0\leq i\leq d-1.
\end{equation}
The intersection numbers are given by
\begin{equation}\label{cub}
b_i=d-i,\;\;\ 0\leq i\leq d-1; \;\;\;\ c_i=i,\;\;\ 1\leq i\leq d.
\end{equation}
Then by using (\ref{l3}) and (\ref{res3}), we obtain
$$R_{\alpha\beta^{(1)}}=\frac{2^d-1}{d2^{d-1}},$$
$$R_{\alpha\beta^{(2)}}=R_{\alpha\beta^{(1)}}+\frac{2^d-1-d}{d(d-1)2^{d-1}}=\frac{2^{d-1}-1}{(d-1)2^{d-2}},\;\;\
$$
\begin{equation}\label{result}
R_{\alpha\beta^{(3)}}=R_{\alpha\beta^{(1)}}+\frac{2^d-1-d-d(d-1)/2}{2^{d-2}d(d-1)(d-2)}=\frac{1}{d(d-1)(d-2)}\{\frac{2^{d}(d^2-2d+2)-3d(d-1)-2}{2^{d-1}}\}.
\end{equation}
The other resistance distances $R_{\alpha\beta^{(i)}}$ for $i\geq
4$, can be calculated similarly.
\subsection{Johnson network}
Let $n\geq 2$ and $d\leq n/2$. The Johnson network $J(n,d)$ has all
$d$-element subsets of $\{1,2,...,n\}$ such that two $d$-element
subsets are adjacent if their intersection has size $d-1$. Two
$d$-element subsets are then at distance $i$ if and only if they
have exactly $d-i$ elements in common. The Johnson network $J(n,d)$
has $N=\frac{n!}{d!(n-d)!}$ vertices, diameter $d$ and the valency
$\kappa=d(n-d)$. Its intersection array is given by
\begin{equation}
b_i=(d-i)(n-d-i), \;\;\;\ 0\leq i\leq d-1; \;\;\ c_i=i^2, \;\;\;\
1\leq i\leq d,
\end{equation}
Then by using (\ref{l3}) and (\ref{res3}), one can obtain
$$R_{\alpha\beta^{(1)}}=\frac{2(n!-d!(n-d)!)}{d(n-d)n!},$$
$$R_{\alpha\beta^{(2)}}=R_{\alpha\beta^{(1)}}+\frac{2[n!-(1+d(n-d))d!(n-d)!]}{n!d(d-1)(n-d)(n-d-1)}=$$
$$\frac{2}{d(d-1)(n-d)(n-d-1)}\{d(n-d)-(n-2)+\frac{d!(n-d)!(n-2-2d(n-d))}{n!}\},$$
$$R_{\alpha\beta^{(3)}}=R_{\alpha\beta^{(2)}}+\frac{2}{\frac{n!}{d!(n-d)!}d(d-1)(d-2)(n-d)(n-d-1)(n-d-2)/4}\{\frac{n!}{d!(n-d)!}-1-$$
$$d(n-d-d(d-1)(d-2)(n-d)(n-d-1)(n-d-2)/4)\}=\frac{2}{d(d-1)(d-2)(n-d)(n-d-1)(n-d-2)}\times$$
$$
\{d^2(n-2d+1)+(3n-2d(n-d)-10)\frac{d(n-d)d!(n-d)!}{n!}+[d^2(n-d)^2-d(n-d)(3n-9)-$$
\begin{equation}\label{resultJ}
4(d-1)(n-d-1)+2(n-2)(n-4)](1-\frac{d!(n-d)!}{n!})\}.
\end{equation}
Again, one can obtain the other resistance distances
$R_{\alpha\beta^{(i)}}$ for $i\geq 4$, similarly.
\section{Conclusion}
Based on the Bose-Mesner algebra corresponding to distance-regular
networks, stratification, spectral techniques and
Christoffel-Darboux identity, a recursive formula for calculating
resistance distance in distance-regular resistor networks was
obtained such that one can evaluate the resistance distances on
these networks only by knowing the corresponding intersection array
, without any need to know the spectrum of the pseudo inverse of the
Laplacian of the networks. As an important result, it was shown that
the resistance distance on a distance-regular network is an
increasing function of the shortest path distance defined on the
network. Although we focused specifically on distance-regular
networks, we hope that the introduced method might then be applied
to other underlying networks of association schemes which are not
distance-regular ones such as finite and infinite square lattice and
underlying networks of the root lattices of type $A_n$ particularly
finite and infinite hexagonal networks ($n=2$) introduced in
\cite{rootlatt} by employing the Krylov-subspace Lanczos algorithm
\cite{js2} iteratively to give three-term recursion relations to the
networks, where these problems are under investigation.
\newpage
 \vspace{1cm}\setcounter{section}{0}
 \setcounter{equation}{0}
 \renewcommand{\theequation}{A-\roman{equation}}
  {\Large{Appendix}}\\
In this appendix, we give the two-point resistances
$R_{\alpha\beta^{(i)}}$, $i=1,2,...,d$ for some important
finite distance-regular networks.\\
\textbf{1. Collinearity graph, gen. octagon $(s,1)$, $s=2,3,4$}\\
$$N=s^4+2s^3+2s^2+2s+1,\;\;\ \{b_0,b_1,b_2,b_3;c_1,c_2,c_3,c_4\}=\{2s,s,s,s;1,1,1,2\},$$
$$\kappa=b_0=2s,\;\ \kappa_2=\frac{\kappa b_1}{c_2}=2s^2,\;\ \kappa_3=\frac{\kappa_2 b_2}{c_3}=2s^3,\;\ \kappa_4=\frac{\kappa_3 b_3}{c_4}=s^4.$$
Then,
$$R_{\alpha\beta^{(1)}}=\frac{s^3+2s^2+2s+2}{s^4+2s^3+2s^2+2s+1},$$
$$R_{\alpha\beta^{(2)}}=R_{\alpha\beta^{(1)}}+\frac{s^2+2s+2}{s^4+2s^3+2s^2+2s+1}=\frac{s^3+3s^2+4s+4}{s^4+2s^3+2s^2+2s+1},$$
$$R_{\alpha\beta^{(3)}}=R_{\alpha\beta^{(2)}}+\frac{s+2}{s^4+2s^3+2s^2+2s+1}=\frac{s^3+3s^2+5s+6}{s^4+2s^3+2s^2+2s+1},$$
$$R_{\alpha\beta^{(4)}}=R_{\alpha\beta^{(3)}}+\frac{1}{s^4+2s^3+2s^2+2s+1}=\frac{s^3+3s^2+5s+7}{s^4+2s^3+2s^2+2s+1}.$$\\
\textbf{2. Incidence graph, $pg(l - 1; l - 1; l - 1),l = 4; 5; 7;
8$}
$$N=2l^2,\;\;\ \{b_0,b_1,b_2,b_3;c_1,c_2,c_3,c_4\}=\{l,l-1,l-1,1;1,1,l-1,l\},$$
$$\kappa=l,\;\ \kappa_2=l(l-1),\;\ \kappa_3=l(l-1),\;\ \kappa_4=l-1.$$
Then,
$$R_{\alpha\beta^{(1)}}=\frac{2l^2-1}{l^3},\;\;\ R_{\alpha\beta^{(2)}}=R_{\alpha\beta^{(1)}}+\frac{2l^2-1-l}{l^3(l-1)}=\frac{2l^2-1}{l^2(l-1)},$$
$$R_{\alpha\beta^{(3)}}=R_{\alpha\beta^{(2)}}+\frac{l+1}{l^3(l-1)}=\frac{2l^3+1}{l^3(l-1)},\;\;\
R_{\alpha\beta^{(4)}}=R_{\alpha\beta^{(3)}}+\frac{1}{l^3}=\frac{2l^2+1}{l^2(l-1)}.$$\\
\textbf{3. Hadamard graph}
$$N=16\gamma,\;\;\ \{b_0,b_1,b_2,b_3;c_1,c_2,c_3,c_4\}=\{4\gamma,4\gamma-1,2\gamma,1;1,2\gamma,4\gamma-1,4\gamma\},\;\ \mathrm{where},\;\ \gamma\in \mathrm{N}$$
$$\kappa=4\gamma,\;\ \kappa_2=2(4\gamma-1),\;\ \kappa_3=4\gamma,\;\ \kappa_4=1.$$
Then,
$$R_{\alpha\beta^{(1)}}=\frac{16\gamma-1}{32\gamma^2},\;\;\ R_{\alpha\beta^{(2)}}=R_{\alpha\beta^{(1)}}+\frac{12\gamma-1}{32\gamma^2(4\gamma-1)}=\frac{8\gamma-1}{4\gamma(4\gamma-1)},$$
$$R_{\alpha\beta^{(3)}}=R_{\alpha\beta^{(2)}}+\frac{4\gamma+1}{32\gamma^2(4\gamma-1)}=\frac{64\gamma^2-4\gamma+1}{32\gamma^2(4\gamma-1)},\;\;\ R_{\alpha\beta^{(4)}}=R_{\alpha\beta^{(3)}}+\frac{1}{32\gamma^2}=\frac{2}{4\gamma-1}.$$
\\
\textbf{4. Distance-regular graphs with classical parameters}. Let
$\Gamma$ denote a distance-regular graph with diameter $d\geq3$. We
say $\Gamma$ has classical parameters $(d, q, \alpha, \beta)$
whenever the intersection numbers are given by
$$c_i=\left[\begin{array}{c}
        i \\
        1
      \end{array}\right] (1+\alpha\left[\begin{array}{c}
        i-1 \\
        1
      \end{array}\right] ),\;\;\ i=1,2,...,d,
$$
\begin{equation}\label{intpar}
b_i=(\left[\begin{array}{c}
        d \\
        1
      \end{array}\right] -\left[\begin{array}{c}
        i \\
        1
      \end{array}\right] )(\beta-\alpha\left[\begin{array}{c}
        i \\
        1
      \end{array}\right] ),\;\;\ i=0,1,...,d-1,
\end{equation}
where, $$\left[\begin{array}{c}
        j \\
        1
      \end{array}\right]:=1+q+q^2+...+q^{j-1}.$$
For instance, the $d$-cube is a distance-regular graph with the
classical parameters $d$, $q=1$, $\alpha=0$ and $\beta=1$. As an
another example of these types of graphs, one can consider the graph
$B_d(p^n)$ which is a type of so-called Dual Polar Graphs
\cite{weng} and is characterized by $\alpha=0,\;\ \beta=p^n,$ and
$$b_i=\frac{q^{i+1}(q^{d-i}-1)}{q-1},\;\ i=0,1,...,d-1$$
\begin{equation}\label{dual}
 c_i=\frac{q^i-1}{q-1},\;\ i=1,2,...,d,
\end{equation}
where, $p$ is a prime number and $n\in \mathrm{N}$. For example,
consider $d=4$, and $q=2$ then we have the graph $B_4(2)$ with the
number of vertices, intersection arrays and valencies
$$N=2295,\;\;\ \{b_0,b_1,b_2,b_3;c_1,c_2,c_3,c_4\}=\{30,28,24,16;1,3,7,15\},$$
$$\kappa=30,\;\ \kappa_2=280,\;\ \kappa_3=960,\;\ \kappa_4=1024.$$
Then,
$$R_{\alpha\beta^{(1)}}=\frac{2294}{34425},\;\;\;\
R_{\alpha\beta^{(2)}}=R_{\alpha\beta^{(1)}}+\frac{566}{34425\times7}=\frac{16623}{240975},\;\;\;\
R_{\alpha\beta^{(3)}}=R_{\alpha\beta^{(2)}}+\frac{62}{240975}=\frac{16685}{240975},$$
$$R_{\alpha\beta^{(4)}}=R_{\alpha\beta^{(3)}}+\frac{2}{34425}=\frac{16699}{240975}.$$\\
\textbf{5. $M_{22}$ graph}
$$N=330,\;\;\ \{b_0,b_1,b_2,b_3;c_1,c_2,c_3,c_4\}=\{7,6,4,4;1,1,1,6\},$$
$$\kappa=7,\;\ \kappa_2=42,\;\ \kappa_3=168,\;\ \kappa_4=112.$$
Then,
$$R_{\alpha\beta^{(1)}}=\frac{47}{165},\;\;\;\
R_{\alpha\beta^{(2)}}=R_{\alpha\beta^{(1)}}+\frac{161}{3465}=\frac{164}{495},\;\;\;\
R_{\alpha\beta^{(3)}}=R_{\alpha\beta^{(2)}}+\frac{1}{99}=\frac{1183}{3465},$$
$$R_{\alpha\beta^{(4)}}=R_{\alpha\beta^{(3)}}+\frac{1}{990}=\frac{113}{330}.$$\\
\textbf{6. Collinearity graph gen. dodecagon $GD(s; 1)$}
$$N=s^6+2(s^5+s^4+s^3+s^2+s)+1,\;\;\ \{b_0,b_1,b_2,b_3,b_4,b_5;c_1,c_2,c_3,c_4,c_5,c_6\}=\{2s, s, s, s, s, s;1, 1, 1, 1, 1, 2\},$$
$$\kappa=2s,\;\ \kappa_2=2s^2,\;\ \kappa_3=2s^3,\;\ \kappa_4=2s^4,\;\ \kappa_5=2s^5,\;\ \kappa_6=s^6.$$
Then,
$$R_{\alpha\beta^{(1)}}=\frac{s^5+2(s^4+s^3+s^2+s+1)}{s^6+2(s^5+s^4+s^3+s^2+s)+1},$$
$$R_{\alpha\beta^{(2)}}=R_{\alpha\beta^{(1)}}+\frac{s^4+2(s^3+s^2+s+1)}{s^6+2(s^5+s^4+s^3+s^2+s)+1}=\frac{s^5+3s^4+4(s^3+s^2+s+1)}{s^6+2(s^5+s^4+s^3+s^2+s)+1},$$
$$R_{\alpha\beta^{(3)}}=R_{\alpha\beta^{(2)}}+\frac{s^3+2(s^2+s+1)}{s^6+2(s^5+s^4+s^3+s^2+s)+1}=\frac{s^5+3s^4+5s^3+6(s^2+s+1)}{s^6+2(s^5+s^4+s^3+s^2+s)+1},$$
$$R_{\alpha\beta^{(4)}}=R_{\alpha\beta^{(3)}}+\frac{s^2+2(s+1)}{s^6+2(s^5+s^4+s^3+s^2+s)+1}=\frac{s^5+3s^4+5s^3+7s^2+8(s+1)}{s^6+2(s^5+s^4+s^3+s^2+s)+1},$$
$$R_{\alpha\beta^{(5)}}=R_{\alpha\beta^{(4)}}+\frac{s+2}{s^6+2(s^5+s^4+s^3+s^2+s)+1}=\frac{s^5+3s^4+5s^3+7s^2+9s+10)}{s^6+2(s^5+s^4+s^3+s^2+s)+1},$$
$$R_{\alpha\beta^{(6)}}=R_{\alpha\beta^{(5)}}+\frac{1}{s^6+2(s^5+s^4+s^3+s^2+s)+1}=\frac{s^5+3s^4+5s^3+7s^2+9s+11}{s^6+2(s^5+s^4+s^3+s^2+s)+1}.$$\\
\textbf{7. A distance-regular graph with}
$$N=2[1+l+\frac{l(l-1)}{c}],\;\;\ \{b_0,b_1,b_2,b_3,b_4;c_1,c_2,c_3,c_4,c_5\}=\{l, l-1, l-c, c, 1;1, c, l-c, l-1, l\},$$
where $l=\gamma(\gamma^2+3\gamma+1)$ and $c=\gamma(\gamma+1),\;\
\gamma\in \mathrm{N}$. Then,
$$\kappa=l,\;\ \kappa_2=\frac{l(l-1)}{c},\;\ \kappa_3=\frac{l(l-1)}{c},\;\ \kappa_4=l,\;\ \kappa_5=1.$$
Then,
$$R_{\alpha\beta^{(1)}}=\frac{c+2lc+2l(l-1)}{l(c+lc+l(l-1))},$$
$$R_{\alpha\beta^{(2)}}=R_{\alpha\beta^{(1)}}+\frac{c+lc+2l(l-1)}{l(l-1)(c+lc+l(l-1))}=\frac{2(c+l)}{c+lc+l(l-1)},$$
$$R_{\alpha\beta^{(3)}}=R_{\alpha\beta^{(2)}}+\frac{c}{l(l-1)(l-c)}=\frac{2l(l-1)(l^2-c^2)+c^2(l+1)+cl(l-1)}{l(l-1)(l-c)(c+cl+l(l-1))},$$
$$R_{\alpha\beta^{(4)}}=R_{\alpha\beta^{(3)}}+\frac{c(l+1)}{l(l-1)(c+lc+l(l-1))}=\frac{2[l^2(2l^2-c^2-l)+cl(l+c)]}{l(l-1)(l-c)(c+lc+l(l-1))},$$
$$R_{\alpha\beta^{(5)}}=R_{\alpha\beta^{(4)}}+\frac{c}{l(c+lc+l(l-1))}=\frac{2l^2(2l^2-c^2-l)+cl(3l+c-1)+c^2}{l(l-1)(l-c)(c+lc+l(l-1))}.$$\\
\textbf{8. Biggs-Smith graph}
$$N=102,\;\;\ \{b_0,b_1,b_2,b_3,b_4,b_5,b_6;c_1,c_2,c_3,c_4,c_5,c_6,c_7\}=\{3, 2, 2, 2, 1, 1, 1 ;1, 1, 1, 1, 1,1,3\},$$
$$\kappa=3,\;\ \kappa_2=6,\;\ \kappa_3=12,\;\ \kappa_4=24,\;\ \kappa_5=24,\;\ \kappa_6=24,\;\ \kappa_7=8.$$
Then,
$$\hspace{-1cm}R_{\alpha\beta^{(1)}}=\frac{101}{153},\;\;\ R_{\alpha\beta^{(2)}}=R_{\alpha\beta^{(1)}}+\frac{49}{153}=\frac{150}{153},\;\;\
R_{\alpha\beta^{(3)}}=R_{\alpha\beta^{(2)}}+\frac{23}{153}=\frac{173}{153},$$
$$R_{\alpha\beta^{(4)}}=R_{\alpha\beta^{(3)}}+\frac{10}{153}=\frac{183}{153},\;\;\ R_{\alpha\beta^{(5)}}=R_{\alpha\beta^{(4)}}+\frac{7}{153}=\frac{190}{153},\;\;\
R_{\alpha\beta^{(6)}}=R_{\alpha\beta^{(5)}}+\frac{4}{153}=\frac{194}{153},$$
$$R_{\alpha\beta^{(7)}}=R_{\alpha\beta^{(6)}}+\frac{1}{153}=\frac{195}{153}.$$\\
\textbf{9. Foster graph}
$$N=90,\;\;\ \{b_0,b_1,b_2,b_3,b_4,b_5,b_6,b_7;c_1,c_2,c_3,c_4,c_5,c_6,c_7,c_8\}=\{3, 2, 2, 2, 2,1, 1, 1 ;1, 1, 1, 1,2,2,2,3\},$$
$$\kappa=3,\;\ \kappa_2=6,\;\ \kappa_3=12,\;\ \kappa_4=24,\;\ \kappa_5=24,\;\ \kappa_6=12,\;\ \kappa_7=6,\;\ \kappa_8=2.$$
Then,
$$\hspace{-2cm}R_{\alpha\beta^{(1)}}=\frac{89}{135},\;\;\;\
R_{\alpha\beta^{(2)}}=R_{\alpha\beta^{(1)}}+\frac{43}{135}=\frac{132}{135},\;\;\;\
R_{\alpha\beta^{(3)}}=R_{\alpha\beta^{(2)}}+\frac{4}{27}=\frac{152}{135},$$
$$\hspace{0.5cm}R_{\alpha\beta^{(4)}}=R_{\alpha\beta^{(3)}}+\frac{17}{270}=\frac{321}{270},\;\;\ R_{\alpha\beta^{(5)}}=R_{\alpha\beta^{(4)}}+\frac{11}{540}=\frac{653}{540},\;\;\
R_{\alpha\beta^{(6)}}=R_{\alpha\beta^{(5)}}+\frac{5}{270}=\frac{663}{540},$$
$$\hspace{-3.5cm}R_{\alpha\beta^{(7)}}=R_{\alpha\beta^{(6)}}+\frac{2}{135}=\frac{671}{540},\;\;\;\ R_{\alpha\beta^{(8)}}=R_{\alpha\beta^{(7)}}+\frac{1}{135}=\frac{675}{540}.$$\\

\newpage
{\bf Figure Captions}

{\bf Figure-1:}(a) Shows the cube or Hamming scheme $H(3,2)$ with
vertex set $\hspace{0.5cm} V=\{(ijk): i,j,k=0,1\}$ and relations
$\hspace{0.5cm} R_0=\{((ijk),(ijk)):(ijk)\in V\},$ $\hspace{0.5cm}
R_1=\{((ijk),(i'jk)),\\((ijk),(ij'k)),((ijk),(ijk')):i\neq
i',j\neq j',k\neq k'\}$, $\hspace{0.5cm}
R_2=\{((ijk),(i'j'k)),((ijk),(i'jk')),\\((ijk),(ij'k')): i\neq
i',j\neq j',k\neq k'\}$ and $R_3=\{((ijk),(i'j'k')):i\neq i',j\neq
j',k\neq k' \}$ respectively. Its non-vanishing intersection
numbers are: $p_{11}^0=3, \;\ p_{11}^2=2,\;\
p_{12}^1=p_{21}^1=2,\;\ p_{12}^3=p_{21}^3=3,\;\
p_{13}^2=p_{31}^2=1,\;\ p_{22}^0=3,\;\ p_{22}^0=3,\;\
p_{23}^1=p_{32}^1=1,  \;\ p_{33}^0=1.$ (b)The vertical dashed
lines denote the four strata of the cube.

{\bf Figure-2:} Shows the octahedron or Johnson scheme $J(4,2)$.

{\bf Figure-3:} Shows edges through $\alpha$ and $\beta$ in a
distance-regular graph.
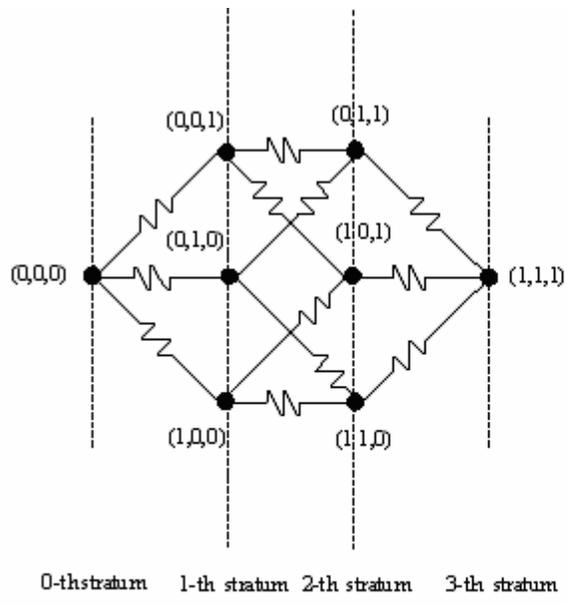
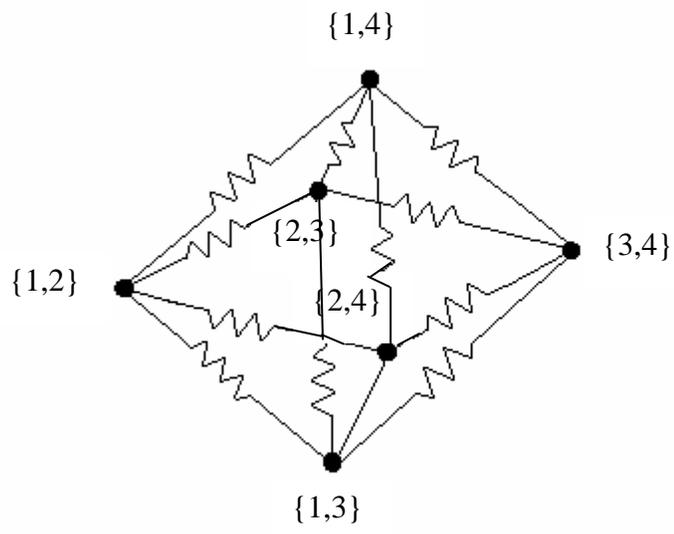
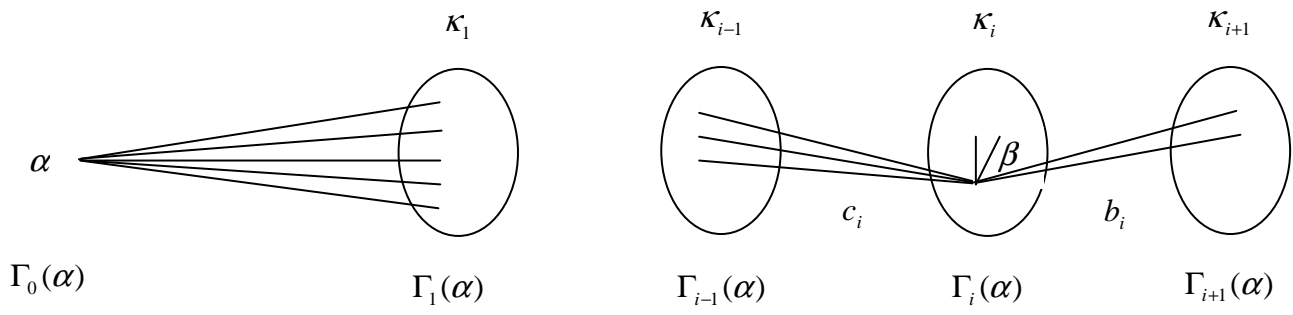
\end{document}